\documentclass[11pt]{article}
\begin{document}
\title{Why I don't like ``pure mathematics''}
\author{\it Volker Runde}
\date{}
\maketitle
I am a pure mathematician, and I enjoy being one. I just don't like the adjective ``pure'' in ``pure mathematics''. Is mathematics that has applications somehow ``impure''? The English mathematician Godfrey Harold Hardy thought so. In his book 
{\it A Mathematician's Apology\/}, he writes:
\begin{quote} \it 
A science is said to be useful of its development tends to accentuate the existing inequalities in the distribution of wealth, or more directly promotes the destruction of human life. 
\end{quote}
His criterion for good mathematics was an entirely aesthetic one:
\begin{quote} \it 
The mathematician's patterns, like the painter's or the poet's must be beautiful; the ideas, like the colours or the words must fit together in a harmonious way. Beauty is the first
test: there is no permanent place in this world for ugly mathematics. 
\end{quote}
I tend to agree with the second quote, but not with the first one.
\par
Hardy's book was written in 1940, when the second world war was raging and the memory of the first one was still fresh. The first world war was the first truly modern war in the sense that science was systematically put to use on the battlefield.
Physicists and chemists helped to develop weapons of unheard of lethal power. After that war, nobody could claim anymore that science was mainly the noble pursuit of knowledge: science had an impact on the real world, sometimes a devastating one, and
scientist could no longer eschew the moral issues involved. By declaring mathematics --- or at least good mathematics --- to be without applications, he absolved mathematics, and thus the mathematical community, from being an accomplice of those
who waged wars and thrived on social injustice.
\par
The problem with this view is simply that it isn't true. Mathematicians live in the real world, and their mathematics interacts with the real world in one way or another. I don't want to say that there is no difference between pure and applied math:
someone who uses mathematics to maximize the time an airline's fleet is actually in the air (thus making money) and not on the ground (thus costing money) is doing applied math whereas someone who proves theorems on the Hochschild cohomology of
Banach algebras (I do that, for instance) is doing pure math. In general, pure mathematics has no {\it immediate\/} impact on the real world (and most of it probably never will), but once we omit that adjective, the line begins to blur.
\par
The fundamental theorem of arithmetic was already known to the ancient Greeks: every positive integer has a prime factorization that is unique up to the order of the factors. A proof is given in Euclid's more than two thousand years old {\it Elements\/}, and there is
little doubt that it was known long before it found its way into that book. For centuries, this theorem was the epitomy of beautiful, but otherwise useless mathematics. This changed in the 1970s with the discovery of the RSA algorithm. It is easy to 
multiply integers on a computer; it is much harder --- even though the fundamental theorem says that it can always be done --- to determine the prime factorization of a given positive integer. This fact can be used to create codes that are extremely hard to 
crack. Without them, e-commerce as it exists today would be impossible. Who would want to key his/her credit card number into an online form if he/she had no guarantee that no eavesdropping crook could get hold of it?
\par
Another mathematical ingredient of the RSA algorithm is Fermat's little theorem (not to be confused with his much more famous last theorem). Pierre de Fermat, a lawyer and civil servant in 17th century France, was doing mathematics in his free time. 
He did it because he enjoyed the intellectual challenge of it, not because it had any connection with his day job. Here is his little theorem: {\it If $p$ is a prime number and $a$ is any integer that does not contain $p$ as a prime factor, then
$p$ divides $a^{p-1} - 1$\/}. This theorem is not obvious, but also not very hard to prove (it probably is on the syllabus of every undergraduate course in number theory). Fermat proved it out of curiosity. Computers, let alone e-commerce, didn't exist
in his days. Nevertheless, it turned out to be useful more than three hundred years after its creator had died.
\par
At the time of Fermat's death, Gottfried Wilhelm von Leibniz was 19 years old. Long after his death, he would be called the last universal genius: he may have been the last person to have a complete grasp of the amassed knowledge of his time.
As a mathematician, he was one of the creators of calculus --- not a small accomplishment ---, and he attempted, but ultimately failed, to build a calculating machine, a forerunner of today's computers. As a philosopher, he gained fame (or notoriety)
through an essay entitled {\it Th\'eodic\'ee\/} --- meaning: God's defense --- in which he tried to reconcile the belief in a loving, almighty God with the apparent existence of human suffering: he argued that we do indeed live in the best of all
possible worlds. Philosophical and theological considerations led him to discover the binary representation of numbers: instead of expressing a number in the decimal system, e.g., $113 = 1 \cdot 10^2 + 1 \cdot 10 + 3 \cdot 10^0$, we can do it equally well 
in the dual system ($113 = 1 \cdot 2^6+1 \cdot 2^5+ 1 \cdot 2^4+0 \cdot 2^3 + 0 \cdot 2^2 + 0 \cdot 2 + 1 \cdot 2^0$). Since numbers in binary representation are much easier to handle on an electronic computer, Leibniz' discovery helped to at least
facilitate the inception of modern information technology.
\par
Almost three hundred years after Leibniz had died, the mathematician Vaughn Jones was working on the problem of classifying subfactors (I won't attempt to explain what a subfactor is; it has nothing to do with multiplying numbers). To accomplish
this classification, he introduced what is now called the Jones index: with each subfactor a certain number is associated. This index displayed a rather strange behavior: it could be infinity or any real number greater than or equal to $4$, but the
values it attained under $4$ had to be of the form $4\cos^2(\pi/n)$
with $n = 3,4, \ldots$. Jones asked himself why. His research led to
the discovery of the Jones polynomial (of course, he didn't call it that way) for which he was awarded the Fields Medal, the
highest honor that can be bestowed upon a (pure) mathematician. This Jones polynomial, in turn, has helped molecular biologists to better understand the ways DNA curls up in a cell's nucleus.
\par
Most of pure mathematics will probably never impact the world outside the mathematical community, but who can be sure in a particular case? In the last twenty five years, the intellectual climate in most ``developed'' countries has become
increasingly unfavorable towards {\it l'art pour l'art\/}. Granting agencies nowadays demand from researchers to explain what the gains of their research are. In principle, there is nothing wrong with that: taxpayers have a right to know what their money is 
used for. The problem is the timeframe: the four examples I gave show that research that was done for nothing but curiosity and the sheer pleasure of exploration, turned out to have applications, sometimes with far reaching consequences. To abandon 
theoretical research just because it doesn't have any foreseeable application in the near future, is a case of cutting off one's nose despite the face.
\par
Pure mathematics isn't pure: neither in the sense that it is removed from the real world, nor in the sense that its practitioners can ultimately avoid the moral questions faced by more applied scientists. It would much better be called
``theoretical mathematics''.
\par
P.S. While Hardy wrote his {\it Apology\/}, other British mathematicians worked on and eventually succeeded in breaking the Enigma code used by the German navy. By all likelihood, their work helped shorten the war by months if not years, thus saving millions of lives on both
sides.
\par
P.P.S. In 1908, Hardy came up with a law that described how the proportions of dominant and recessive genetic traits are propagated in large populations. He didn't think much of it, but it has turned out to be of major importance in the study of blood 
group distributions. 
\vfill
\begin{tabbing}
{\it Address\/}: \= Department of Mathematical and Statistical Sciences \\
\> University of Alberta \\
\> Edmonton, Alberta \\
\> Canada T6G 2G1 \\[\medskipamount]
{\it E-mail\/}: \> {\tt vrunde@ualberta.ca} \\[\medskipamount]
{\it URL\/}: \> {\tt http://www.math.ualberta.ca/$^\sim$runde/}
\end{tabbing} 
\end{document}